\documentclass[a4paper,12pt]{amsart}

\usepackage{amscd, amssymb, amsmath}

\input{prepictex}
\input{pictex}
\input{postpictex}

\def\is{\text{intrinsically synchronising}}

\newtheorem{thm}{Theorem}[section]
\newtheorem{cor}[thm]{Corollary}
\newtheorem{lem}[thm]{Lemma}
\newtheorem{prop}[thm]{Proposition}

\theoremstyle{definition}
\newtheorem{dfn}[thm]{Definition}
\newtheorem{dfns}[thm]{Definitions}

\theoremstyle{remark}
\newtheorem{rmk}[thm]{Remark}
\newtheorem{rmks}[thm]{Remarks}
\newtheorem{example}[thm]{Example}
\newtheorem{examples}[thm]{Examples}

\title{Reducibility of Covers of AFT shifts}

\author{Teresa Bates}
\address{Teresa Bates\\ School of Mathematics \& Statistics \\ University
of New South Wales \\ UNSW Sydney NSW 2052 \\ AUSTRALIA} \email{teresa@unsw.edu.au}

\author{S\o ren Eilers}
\address{S\o ren Eilers\\ Department of Mathematical Sciences \\
Universitetsparken 5 \\
2100 Copenhagen \O \\ DENMARK} \email{eilers@math.ku.dk}

\author{David Pask}
\address{David Pask \\ School of Mathematics \&
Applied Statistics  \\
University of Wollongong\\
NSW  2522\\
AUSTRALIA} \email{dpask@uow.edu.au}



\thanks{This research was supported by the UNSW Faculty Research Grants, Australian Research Council.}

\date{\today}

\begin{document}

\begin{abstract}
In this paper we show that the reducibility structure of several covers of
sofic shifts is a flow invariant.  In addition, we prove
that for an irreducible subshift of almost finite type the left Krieger cover
and the past set cover are reducible.  We provide an example which shows that
there are non almost finite type shifts which have reducible left Krieger covers.
As an application we show that the Matsumoto algebra of an irreducible, strictly sofic
shift of almost finite type is not simple.
\end{abstract}

\maketitle

\section{Introduction}

The classification of shift spaces of finite type up to flow equivalence,
initiated by Parry and Sullivan and Bowen and Franks (see \cite{parsu,bf}) and completed by Boyle
and Huang \cite{bh} is, at present, not generalized to any but a few sporadic
classes of shift spaces. This state of affairs is certainly related to
the scarcity of useful invariants presently known, even in the case when
the shift space is irreducible. The present paper follows a general
strategy, also pursued in, for example \cite{ce1,ce2,m6}, of trying to extract
such invariants from the study of $C^*$-algebras associated to shift spaces:
since the $C^*$-algebra is a flow invariant, anything that derives from it must
also be one. However, we can (and will) suppress $C^*$-algebra theory from the
presentation in this paper.

The class of sofic shifts coined by Weiss \cite{wei} and subsequently studied
intensively by several authors (see, for example, \cite{bfp,blan,ffj,fisc,fo,jm,kri,kri1,nasu,thoms3})
captures, in a sense, the next level of complexity after the shifts of finite type.  The purpose of this paper
is to investigate the reducibility structure of various left-resolving
presentations of a sofic shift space and prove that this structure is a flow invariant.
The most commonly used presentations
that we refer to are the left Krieger cover, the left Fischer cover and the
past set cover.   The precise relationships between these covers are currently
unclear to us, but will be explored further in \cite{rune}.

In this paper we work towards providing conditions on an irreducible sofic shift
which guarantee that the left Krieger (respectively past set) cover is reducible.
Theorem \ref{thm:aftresult} and Theorem \ref{thm:aftresult2} do not characterise the
irreducible sofic shifts having this property as Examples \ref{ex:bfg} and
Example \ref{ex:jonathanstrikesback} demonstrate.

In Section \ref{sec:bg} we recall the background theory of sofic shift spaces, left Krieger
covers and shifts of almost finite type.  The definition of the left Krieger cover we take
is the one given in \cite{kri,kit}, but we note that in some places
in the literature (for example, \cite{c,cm})  the left Krieger cover is defined
to be the cover $(E^\infty_f, {\mathcal L})$ discussed in Remarks \ref{rmk:othercover} (i).
 Our intention here is two-fold:  first we wish to make
the present paper self-contained and second we hope to clarify some inconsistencies in the
terminology found in the literature which have caused us confusion.

In Section \ref{sec:flowinv} we show that the reducibility structure of the left Krieger
cover is invariant under the operations that generate flow equivalence for sofic shifts.
Our main contribution here is to show in Proposition \ref{prop:symexpand}
that the irreducibility structure
of the left Krieger cover of a general subshift is invariant under symbol expansion.  
Krieger has shown in \cite{kri} that the left Krieger cover of a sofic shift
is a conjugacy invariant, and this, in combination with our result, gives invariance of the
reducibility structure of the left Krieger cover of a sofic shift under flow equivalence.

In Section \ref{sec:aftresult} we prove our main theorem.  In Theorem \ref{thm:aftresult}
we show that the left Krieger cover of an irreducible strictly sofic shift of almost finite
type is reducible.   As an application of Theorem \ref{thm:aftresult}
we prove that the Matsumoto algebras (see \cite{m},\cite{cm})
of  irreducible almost finite type shift spaces are not simple in Corollary \ref{cor:notsomatso}.
We show in Example \ref{ex:jonathanstrikesback} that the converse of Theorem \ref{thm:aftresult}
does not hold:  there are examples of irreducible strictly sofic shifts which are not of almost finite
type but which have reducible left Krieger covers.  We also briefly discuss another left-resolving
cover $(E^f_\infty,{\mathcal L})$ of a two-sided sofic shift space which is constructed using the
infinite pasts of finite words.  We note in Remark \ref{rmk:infinitepastfinitefuture} that it is possible to prove a
corresponding version of  Theorem \ref{thm:aftresult} for this cover.

In Section \ref{sec:aftresult2} we note that the construction of the above two covers does
not make sense for one-sided shift spaces.  We consider the past set cover of a (one- or two-sided)
sofic shift space and prove in Theorem \ref{thm:aftresult2} that this cover is reducible for
irreducible strictly sofic shifts of almost finite type.  We note in Remarks \ref{rmk:othercover} that
there are examples of irreducible two-sided strictly sofic shifts for which the left Krieger cover and
the past set cover do not coincide.  Indeed, it is possible for the left Krieger cover
to be irreducible when the past set cover is reducible. In Remarks \ref{rmk:othercover} (i) we consider
a cover of a (one- or two-sided) sofic shift which is constructed by considering the finite pasts of
right-infinite rays in the shift space and give a corresponding version of Theorem \ref{thm:aftresult}
for this cover.  We note that this cover coincides with the left Krieger
cover for two-sided sofic shift spaces.  As an application we are able to give a different proof
of \cite[Proposition 2.7 c)]{thoms4} which states that strictly sofic $\beta$-shifts are not of almost
finite type.

\subsection*{Acknowledgements}
The authors wish to thank Klaus Thomsen and Ian Putnam for some helpful discussions at the Banff International
Research Station. We also wish to thank the Fields Institute in Toronto, Canada where some of the research
for this paper was carried out.
\section{Background} \label{sec:bg}

\subsection*{Sofic shift spaces}


For an excellent treatment of shift spaces we refer the reader to \cite{lm}.  We present
a few basic definitions here in order to make the present paper self-contained.

Let ${\mathcal A}$ be a finite alphabet, and let $\mathcal{A}^*$ denote the free monoid generated by 
$\mathcal{A}$. A \emph{language} is a submonoid of $\mathcal{A}^*$ for some $\mathcal{A}$. The {\em full shift} over
${\mathcal A}$ consists of the space ${\mathcal A}^{\bf Z}$  endowed with the
product topology together with the shift map $\sigma: {\mathcal A}^{\bf Z} \to
{\mathcal A}^{\bf Z}$ such that for $x \in {\mathcal A}^{\bf Z}$, $\sigma(x)_i = x_{i+1}$ for all $i \in {\bf
Z}$.  A {\em subshift} $\textsf{X}$ is a closed $\sigma$-invariant
subset of ${\mathcal A}^{\bf Z}$.  The language of the shift space
$\textsf{X}$ is denoted by ${\mathcal B}(\textsf{X})$ and is the collection of all
{\em words} or {\em blocks} which appear in the bi-infinite sequences
of $\textsf{X}$; the empty word is denoted $\epsilon$.

For each subshift there is a
collection ${\mathcal F}$ of blocks which are not permitted to occur
in the sequences of $\textsf{X}$.  This collection of forbidden
blocks uniquely describes the subshift (cf. \cite[Proposition 1.3.4]{lm}).
If the set of forbidden blocks for $\textsf{X}$ can be chosen to be finite, then
$\textsf{X}$ is called a {\em shift of finite type}  or SFT.

Let $\textsf{X}_1$ and $\textsf{X}_2$ be two subshifts over possibly
different alphabets.   A {\em factor map} $\pi: \textsf{X}_1 \to
\textsf{X}_2$ is a continuous surjective function $\pi: \textsf{X}_1
\to \textsf{X}_2$ which commutes with the shift maps.  According to Weiss \cite{wei}
a shift space is {\em sofic} if it is the image of an SFT
under a factor map.

A directed graph $E$ consists of a quadruple $( E^0 , E^1
, r , s )$ where $E^0$ and $E^1$ are countable sets of vertices
and edges respectively and $r, s : E^1 \to E^0$ are maps giving
the direction of each edge. A path $\lambda = e_1 \ldots e_n$ is a
sequence of edges $e_i \in E^1$ such that $r ( e_i ) = s ( e_{i+1}
)$ for $i=1 , \ldots , n-1$.  We denote the collection of all finite paths in
$E$ by $E^*$ and extend the range and source maps to $E^*$ in the natural way.  
A {\em circuit} in a directed graph $E$ is a finite
path $\lambda \in E^*$ satisfying $s(\lambda) = r(\lambda)$.
A directed graph $E$ is {\em irreducible} if for every pair $(u,v)$ of
vertices there is a path $\lambda \in E^*$ from $u$ to $v$.  Otherwise, the graph is
said to be {\em reducible}.

A directed graph is \emph{essential} if every vertex receives and emits an edge.
A directed graph is \emph{row-finite} if every vertex emits finitely many edges.
We shall work exclusively with essential row-finite graphs.
The edge shift $( \textsf{X}_E , \sigma_E )$ associated
to an essential directed graph $E$ is given by:
\[
\textsf{X}_E = \{ x \in ( E^1 )^{\bf Z} : s ( x_{i+1} ) = r (
x_{i} ) \text{ for all } i \in {\bf Z} \} \text{ and } ( \sigma_E
x )_i = x_{i+1} \text{ for } i \in {\bf Z} .
\]

\noindent
Evidently ${\mathcal B}(\textsf{X}_E) \backslash \{ \epsilon \}= E^* \backslash E^0$.
The following definition is adapted from
\cite[Definition 3.1.1]{lm}:

\begin{dfn}
A {\it labelled graph} $( E , {\mathcal L} )$ over an alphabet ${\mathcal
A}$ consists of a directed graph $E$ together with a surjective labelling map
$\mathcal{L}  : E^1 \to \mathcal{A}$.  We say that the labelled graph $(E,{\mathcal L})$
is {\em essential} if the directed graph $E$ is essential.
\end{dfn}

\noindent
Observe that a directed graph $E$ is a labelled graph $(E,{\mathcal L}_t)$ over the alphabet $E^1$
where ${\mathcal L}_t : E^1 \to E^1$ is the identity map. 

\begin{dfns}  Let  $E,F$ be directed graphs.  A map $\phi=(\phi^0,\phi^1): E \to F$ is an {\em isomorphism
of directed graphs} if $\phi^0:  E^0 \to F^0$ and $\phi^1 : E^1 \to F^1$ are bijections such that for 
all $e \in E^1$ $\phi^0(s(e)) = s(\phi^1(e))$ and $\phi^0(r(e)) = r(\phi^1(e))$.

Let $(E , \mathcal{L}  )$ and $(F , \mathcal{L} ' )$  be labelled graphs over the same alphabet. A graph 
isomorphism $\phi : E \to F$ is a {\em labelled graph isomorphism} if 
$\mathcal{L} ' ( \phi (e) ) = \mathcal{L}  (e)$ for all $e \in E^1$
\end{dfns}

\noindent
Given an essential labelled graph $(E, \mathcal{L})$ over the alphabet ${\mathcal A}$
we may define a subshift $(\textsf{X}_{(E,\mathcal{L})} , \sigma )$ of $\mathcal{A}^{\bf Z}$ by
\[
\textsf{X}_{(E,\mathcal{L} )} = \{ y \in \mathcal{A}^{\bf Z} : \text{ there exists
 } x \in \textsf{X}_E \mbox{ such that } y_i = \mathcal{L} ( x_i ) \text{ for all } i \in {\bf Z}
\} ,
\]

\noindent where $\sigma$ is the shift map inherited from ${\mathcal A}^{\mathbf Z}$.
A {\em representative} of a word $w \in {\mathcal B}(\textsf{X}_{(E,{\mathcal L})})$ is
a path $\lambda = e_1 \ldots e_n \in E^*$ such that
${\mathcal L}( \lambda ) := \mathcal{L} ( e_1 ) \ldots \mathcal{L} ( e_n ) = w$.
The labelled graph $(E, \mathcal{L} )$ is said to be a {\em presentation} of the shift space
$\textsf{X} = \textsf{X}_{(E , \mathcal{L} )}$. As shown in Example \ref{ex:evenshift} a shift
space may have many different presentations.  We shall consider five such
presentations in this paper.

Fischer proved in \cite[Theorem 1]{fisc} that sofic shifts are precisely those shift
spaces that can be presented by labelled
graphs with finite edge and vertex sets.  Since, up to conjugacy,
shifts of finite type are precisely the edge shifts associated to directed graphs with finite edge
and vertex sets (see \cite[Proposition 2.3.9]{lm}), shifts of finite type are sofic shifts as well.  We say that a sofic
shift is {\em strictly sofic} if it is not an SFT.

\begin{dfn}
A shift space $\textsf{X}$ is said to be {\em irreducible} if for every $u,w \in {\mathcal B}(\textsf{X})$
there is a $v \in {\mathcal B}(\textsf{X})$ such that $uvw \in {\mathcal B}(\textsf{X})$.
\end{dfn}

\begin{rmk} \label{rem:soficirred}
It is straightforward to see that a sofic shift is irreducible if and only
if it can be presented by an irreducible labelled graph i.e.\ a labelled
graph $(E,{\mathcal L})$ with $E$ irreducible (see \cite[Section 3.1]{lm} for details).
\end{rmk}

\noindent In this paper we will be primarily
concerned with irreducible strictly sofic shifts.

Let $(E,{\mathcal L})$ be a labelled graph (sometimes called a Shannon graph after \cite{sha})
with finite edge and vertex sets
which presents the sofic shift $\textsf{X}_{(E,{\mathcal L})}$.
By removing the labels from the edges of $(E,{\mathcal L})$ we obtain a presentation
of an SFT $\textsf{X}_E$ and a factor map $\pi_{\mathcal L}: \textsf{X}_E \to \textsf{X}_{(E,{\mathcal L})}$
induced by the $1$-block map ${\mathcal L} : E^1 \to {\mathcal A}$.  The subshift $\textsf{X}_E$ is called a {\em cover} of
the sofic shift $\textsf{X}_{(E,{\mathcal L})}$.

\begin{dfn}
Let $\textsf{X}$ be a shift space.  A word $w \in {\mathcal
B}(\textsf{X})$ is {\em intrinsically synchronising} if whenever $uw, wv \in
{\mathcal B}(\textsf{X})$ we have $uwv \in {\mathcal B}(\textsf{X})$.
\end{dfn}

\begin{rmks} \label{rmk:isremarks}
There  is confusion in the literature as to the correct terminology
for the above concept.  Intrinsically synchronising is called {\em
magic} in \cite{bkm}, {\em finitary} in \cite{kri}  and {\em synchronising} in \cite{trow}.
However these words have different meanings in other places in the
literature (see, for example, \cite{lm}). We shall follow
the terminology used in \cite{lm}.

If $m \in \mathcal{B} ( \textsf{X} )$ is \is \  then for any $u , v \in \mathcal{B} ( \textsf{X} )$
with $umv \in \mathcal{B} ( \textsf{X} )$ the word $umv$ is \is.
\end{rmks}

\begin{dfn}
A labelled graph $(E,{\mathcal L})$ is {\em left-resolving} if for every $v \in E^0$, all edges
ending at $v$ carry different labels.  That is, ${\mathcal L} : r^{-1}(v) \to {\mathcal A}$ is injective
for all $v \in E^0$. A labelled graph is {\em right-resolving} if all
edges leaving each vertex carry different labels.
\end{dfn}

\begin{dfn} \label{def:closingdef}
A labelled graph $(E,{\mathcal L})$ is {\em right-closing with delay} $D$ if every pair of paths
in $E$ of length $D+1$ which start at the same vertex and represent the same word
must have the same initial edge.  That is, whenever paths $\mu,\nu \in E^*$ of length $D+1$
satisfy $s(\mu) = s(\nu)$ and ${\mathcal L}(\mu) = {\mathcal L}(\nu)$ we must have $\mu_1 = \nu_1$.  A labelled
graph is {\em right-closing} if it is right-closing with some delay $D \ge 0$.  The concept of
{\em left-closing} for a labelled graph is similarly defined.
\end{dfn}

\begin{rmk} \label{rem:delaypath}
Let $(E , \mathcal{L} )$ be a labelled graph which is right-closing with delay $D$ and let $\mu , \nu$ be paths
in $E$ of length $D+n+1$ such that $s ( \mu ) = s ( \nu )$ and $\mathcal{L} ( \mu ) = \mathcal{L} ( \nu )$.
Repeated applications of  Definition \ref{def:closingdef} show that
$\mu_1 \ldots \mu_n = \nu_1 \ldots \nu_n$.
\end{rmk}

\subsection*{Left Krieger covers}

Let $\textsf{X}$ be a shift space.

\begin{dfns} (See \cite[Sections I and III]{jm},\cite[Exercise 3.2.8]{lm})\label{def:rightrays}
We write each $x \in \textsf{X}$ as
$x = x^- x^+$ where $x^- = \dots x_{-2}x_{-1}$ is called the {\em left-ray of $x$}  and
$x^+ = x_0 x_1 x_2 \dots$ is called the {\em right-ray of $x$}. We
denote by $\textsf{X}^+$ the set of all right-rays of elements of $\textsf{X}$ and by
$\textsf{X}^-$ the set of all left-rays of elements of $\textsf{X}$.
For $x^+ \in \textsf{X}^+$, the {\em predecessor set} of $x^+$ is the set of all left-rays
$y^-$ which may precede $x^+$; that is
\[
P_\infty ( x^+ ) = \{ y^- \in \textsf{X}^- \, : \, y^-  x^+ \in \textsf{X} \}.
\]
\end{dfns}

\noindent
Note that $\textsf{X}$ is a sofic shift if and only if the number of predecessor sets
$P_\infty ( x^+ )$ is finite, (see \cite[\S 2]{kri}).

\begin{dfn} \label{infinitekrieger}
The {\em left Krieger cover} of a shift space $\textsf{X}$ is the labelled graph
$(E_K,{\mathcal L}_K)$ whose vertices are the predecessor
sets of the elements of $\textsf{X}^+$. There is an edge labelled
$a \in {\mathcal A}$ from $P_\infty(x^+)$ to $P_\infty(y^+)$ if and only if  $ay^+ \in \textsf{X}^+$ and
$P_\infty ( x^+ ) = P_\infty( a y^+ )$.
\end{dfn}

\begin{rmks} \label{rem:lkc}
\begin{itemize}
\item[(i)] If $w \in \mathcal{B} (\textsf{X})$ then for all $x^+ \in \textsf{X}^+$ with $w x^+ \in \textsf{X}^+$ there is
a path labelled $w$ in $( E_K , \mathcal{L}_K )$ beginning at $P_\infty ( w x^+ )$ and
ending at $P_\infty ( x^+ )$. In general, a word $w \in \mathcal{B} ( \textsf{X} )$ can have representatives
which start at several different vertices in the left Krieger cover (see Example \ref{ex:evenshift}).
\item[(ii)] The left Krieger cover of a shift space $\textsf{X}$ is evidently a left-resolving presentation of $\textsf{X}$.
The left Krieger cover is also known as the past state chain, and was originally defined only for sofic shifts by Krieger, but
the definition applies in general.
\item[(iii)] One may similarly define a right Krieger cover (future state chain) which is right-resolving. All of the results about
shift spaces stated in this paper have right-resolving analogues for which the proofs are similar to those given.
\end{itemize}
\end{rmks}

\begin{dfns}(See also \cite[Section 3]{jm})
We say that a ray $x^+ \in \textsf{X}^+$ is {\em intrinsically synchronising} if it
contains an intrinsically synchronising block.  We say that a predecessor set $P_\infty(x^+)$ is
{\em intrinsically synchronising}
if there is an intrinsically synchronising ray $y^+ \in \textsf{X}^+$ with
$P_\infty(x^+) = P_\infty(y^+)$.

We define the {\em left Fischer cover} of a sofic
shift $\textsf{X}$ to be the minimal left-resolving labelled graph $(E_F  , \mathcal{L}_F )$
which presents $\textsf{X}$.
\end{dfns}

\begin{rmks} \label{rem:lfc}
\begin{itemize}
\item[(i)] The right Fischer cover of a sofic shift was originally described in \cite{fisc}
as the minimal right-resolving presentation of a sofic shift.
Note that by a suitably adapted version of \cite[Theorem 3.3.18]{lm} the left Fischer cover
of a sofic shift is well-defined, up to labelled graph isomorphism.

\item[(ii)]The left Fischer cover $(E_F,{\mathcal L}_F)$ of an irreducible sofic shift $\textsf{X}$ may be identified with the
minimal irreducible subgraph of the left Krieger cover of $\textsf{X}$ in which the
vertices are precisely the intrinsically synchronising predecessor sets for the shift
$\textsf{X}$. This result can be traced back to \cite[Lemma 2.7]{kri} however
the proof of Lemma \ref{lem:blmagogo} may be adapted to show this.
\end{itemize}
\end{rmks}

\begin{example} \label{ex:evenshift}
The even shift is an irreducible strictly sofic shift which consists
of the collection of all bi-infinite sequences of $0$'s and $1$'s
such that the number of $0$'s which can occur between two $1$'s is an even number (including zero).

Following \cite[Examples 3.3 (iii)]{bp2}, we see that
the following graphs \newline \goodbreak
$
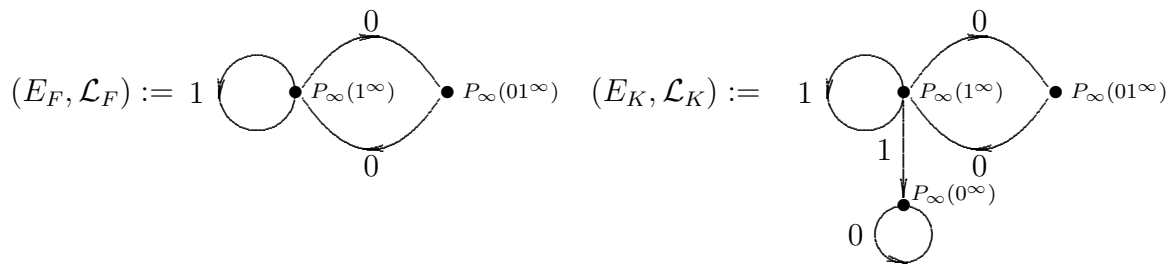
\beginpicture
\setcoordinatesystem units <1cm,0.75cm>

\setplotarea x from 2 to 5, y from -0.2 to 1.6

\put{$(E_F , \mathcal{L} _F ) :=$}[l] at -0.75 0

\put{$\bullet$} at 3 0

\put{$\bullet$} at 5 0

\put{{\tiny $P_\infty(1^\infty)$}}[l] at 3.2 0

\put{{\tiny $P_\infty(01^\infty)$}}[l] at 5.2 0

\put{$1$}[l] at 1.6 0

\put{$0$}[b] at 4 1.1

\put{$0$}[t] at 4 -1.1


\setquadratic

\plot 3.1 0.1  4 1 4.9 0.1 /

\plot 3.1 -0.1 4 -1  4.9 -0.1 /

\circulararc 360 degrees from 3 0 center at 2.5 0

\arrow <0.25cm> [0.1,0.3] from 2.015 0.1 to 2 -0.1 

\arrow <0.25cm> [0.1,0.3] from 4.1 -0.985 to 3.9 -1

\arrow <0.25cm> [0.1,0.3] from 3.9 0.985 to 4.1 1

\endpicture
$ \ \ $
\beginpicture
\setcoordinatesystem units <1cm,0.75cm>

\setplotarea x from 0 to 4, y from -3.5 to 1.6




\put{$(E_K , \mathcal{L} _K ) :=$} at 0 0

\put{$\bullet$} at 3 0

\put{$\bullet$} at 5 0

\put{$\bullet$} at 3 -2

\put{$1$}[l] at 1.6 0

\put{$0$}[b] at 4 1.1

\put{$0$}[t] at 4 -1.1

\put{$1$}[r] at 2.85 -1

\put{$0$}[r] at 2.47 -2.5

\put{{\tiny $P_\infty(1^\infty)$}}[l] at 3.2 0

\put{{\tiny $P_\infty(01^\infty)$}}[l] at 5.2 0

\put{{\tiny $P_\infty(0^\infty)$}}[l] at 3.1 -1.8

\setquadratic

\plot 3.1 0.1  4 1 4.9 0.1 /

\plot 3.1 -0.1 4 -1  4.9 -0.1 /

\circulararc 360 degrees from 3 0 center at 2.5 0

\circulararc 360 degrees from 3 -2 center at 3 -2.5

\arrow <0.25cm> [0.1,0.3] from 2.015 0.1 to 2 -0.1 

\arrow <0.25cm> [0.1,0.3] from 4.1 -0.985 to 3.9 -1

\arrow <0.25cm> [0.1,0.3] from 3.9 0.985 to 4.1 1

\arrow <0.25cm> [0.1,0.3] from 3 -0.15 to 3 -1.85

\arrow <0.25cm> [0.1,0.3] from 2.95 -2.98 to 3.04 -3.01


\endpicture
$\newline
are the left Fischer and left Krieger covers of the even shift, respectively.  The predecessor sets $P_\infty(1^\infty)$ and
$P_\infty(01^\infty)$ are intrinsically synchronising since the rays $1^\infty$ and $01^\infty$ both contain the intrinsically
synchronising word $1$.  Note also that the left Fischer cover shares the vertices $P_\infty(1^\infty)$ and $P_\infty(01^\infty)$
with the left Krieger cover.   These are the only vertices in the left Krieger cover corresponding to intrinsically
synchronising predecessor sets.
\end{example}

\subsection*{AFT shifts}

The original definition of an almost finite type (AFT) shift is due to Marcus
\cite[Definition 4]{marc}.

\begin{dfn} \label{def:AFTdef}
The shift space $\textsf{S}$ is said to be of \emph{almost finite type} if
there is an irreducible subshift of finite type $\textsf{X}$ and
a factor map $\pi: \textsf{X} \to \textsf{S}$ that is one-to-one on a non-trivial open set.
\end{dfn}

In \cite{bfp} a hierarchy of sofic shifts is given, in terms of a certain degree.
By \cite[Proposition 3]{bfp} the irreducible strictly sofic shifts of
degree $1$ are precisely the AFT shifts.  Many authors (see \cite{marc,bkm,will,jm}, for example) have studied AFT shifts and, as a result, we now have the following list of
equivalent conditions on a strictly sofic shift.

\begin{thm} \label{thm:aftconditions}
Let $\textsf{S}$ be a strictly sofic shift.  The following are equivalent
\begin{enumerate}
\item  The shift $\textsf{S}$ is AFT.
\item  The left Fischer cover of $\textsf{S}$ is right-closing.
\item $\textsf{S}$ has a minimal cover (i.e.\ an SFT $\textsf{X}$ and a factor map $\pi: \textsf{X} \to \textsf{S}$ such that any other factor map $\phi: \textsf{B} \to \textsf{S}$ (from an SFT $\textsf{B}$ onto $\textsf{S}$) must factor through $\pi$).
\item  The right and left Fischer covers of $\textsf{S}$ are conjugate
as SFTs.
\end{enumerate}
\end{thm}

\begin{example} The even shift described in Example \ref{ex:evenshift} is an example of an AFT shift since its
left Fischer cover is right-closing.  However, note that its left Krieger cover is reducible and it is
examples of this type which inspire the following sections.
\end{example}

\section{The reducibility structure of left Krieger cover of a shift is a flow invariant} \label{sec:flowinv}

The remarkable result of Parry and Sullivan in \cite{parsu} shows that flow equivalence between shift spaces
is generated by conjugacy and a certain operation that is nowadays referred to as \emph{symbol expansion}. 
In this section we begin by describing the proper communication graph $PC(E)$ of a directed graph $E$, which describes
the irreducible components of the associated shift of finite type $( \textsf{X}_E , \sigma)$. Our main result, Theorem \ref{thm:flowinv}
shows that the proper communication graph of the left Krieger cover of a shift space is a flow invariant -- hence
the reducibility of the left Krieger cover of a shift is a flow invariant.

The following definition is based on \cite[Section 4.4]{lm} (see also \cite[\S 1]{sen}):

\begin{dfn}  Let $v,w$ be vertices in the directed graph $E$.  We say that $v$ {\em is connected
to} $w$ (written $v \geq w$) if there is a path $\mu \in E^*$ with $s(\mu) = v$ and $r(\mu) =  w$.  If
$v,w \in E^0$ satisfy $v \ge w$ and $w \ge v$ then we say that $v$ {\em communicates with} $w$.
\end{dfn}

\begin{rmk}
Each vertex $v$ communicates with itself via the empty path.  It is then straightforward to check that
communication is an equivalence relation. 
\end{rmk}

\noindent
The vertices
of $E$ are partitioned into {\em communicating classes}, that is, maximal sets
of vertices such that each vertex communicates with all the other vertices in the
class.   For each communicating class $C_i$ of $E$, the subgraph $E_i$ of $E$ whose vertices
are the elements of $C_i$ and whose edges are the elements $e$ of $E^1$ with $s(e), r(e)
\in C_i$ is called an {\em irreducible component} of $E$.  If the number of irreducible components
of $E$ is $1$ then $E$ is irreducible.  If an edge $e \in E^1$ has $s(e) \in C_i$ and $r(e) \in C_j$
where $i \ne j$ then we call $e$ a {\em transitional edge for $E$}.
We shall consider the following subset of the communication relation which is not, in general,
an equivalence relation.

\begin{dfn}  We say that vertices $v$ and $w$ {\em properly communicate} if there
are paths $\mu,\nu \in E^*$ of length greater than or equal to one with
$s(\mu) = v, s(\nu) = w$, $r(\mu) = w$ and $r(\nu) = v$.
\end{dfn}

\begin{rmk}
The proper communication relation is not, in general, reflexive.  A vertex
$v$ can only properly communicate with itself if it lies on some circuit $\lambda \in E^*$.
\end{rmk}

\noindent
The proper communication relation may be used to define maximal disjoint subsets
$PC_i$ of vertices such that for all $v,w \in PC_i$ we have that $v$ properly communicates
with $w$.  We call these sets {\em proper communication sets of vertices} and use this
relation to define a directed graph which summarises the structure of the graph $E$ in the following
manner.

\begin{dfn}  Let $E$ be a directed graph.  The {\em proper communication graph} $PC(E)$ is a directed
graph constructed using the following data.  For each proper communication set $PC_i$ of $E$, we draw a vertex
$PC_i$ so that
\[
PC(E)^0 = \{ PC_i \,:\, PC_i \mbox{ is a proper communication set of vertices in } E \}.
\]
For each $i \ne j$ we draw an edge $e_{ij}$ from $PC_i$ to $PC_j$ if there are vertices $v \in PC_i$ and $w \in PC_j$
such that $v \geq w$ so that
\[
PC(E)^1 = \{ e_{ij} \,:\, v \ge w \mbox{ for some } v \in PC_i \mbox{ and } w \in PC_j \}
\]
and we have  $s(e_{ij}) = PC_i$ and $r(e_{ij}) = PC_j$.
\end{dfn}

\begin{rmks}
\begin{enumerate}
\item Note that if $i,j$ are such that there are vertices $v \in PC_i$ and $w \in PC_j$ such that $v \geq w$, then
this relation holds for every pair of vertices $x \in PC_i$ and $y \in PC_j$ since all vertices in a proper
communication set communicate with each other.

\item The proper communication graph $PC(E)$ contains no circuits. Suppose there is a circuit  $\alpha \in PC(E)^*$. Then our construction of $PC(E)$ ensures that $\alpha$ passes through at least two distinct vertices $PC_i$
and $PC_j$.  But then there are vertices $v \in PC_i$ and $w \in PC_j$ such that $v \ge w$ and vertices $x \in PC_j$
and $y \in PC_i$ with $x \geq y$.  Our preceding remark then tells us that we must have $w \ge v$ so that $v$ and
$w$ properly communicate, and so are in the same proper communication set, contradicting the maximality of
$PC_i$ and $PC_j$.

\item  The proper communication graph of row-finite graph $E$ may be used to describe the lattice of gauge-invariant
ideals of the graph $C^*$-algebra $C^*(E)$ described in \cite{kprr,bprsz}, for example.  If $(E_K ,{\mathcal L}_K )$ is the
left Krieger cover of a sofic shift $\textsf{X}$, then by \cite{c} this implies that the proper communication graph
of $E_K$ describes the lattice of gauge-invariant ideals of the Matsumoto algebra ${\mathcal O}_{\textsf{X}}$ described in \cite{m,cm}, for example.
\end{enumerate}
\end{rmks}

\begin{examples}
\begin{enumerate}
\item If the directed graph $E$ contains no circuits, then the proper communication graph $PC(E)$ is
the empty graph.
\item If the directed graph $E$ contains at least one circuit and is irreducible, then the proper
communication graph $PC(E)$ will consist of one vertex and no edges.
\item We draw the proper communication graph $PC(E)$ for the following directed graph
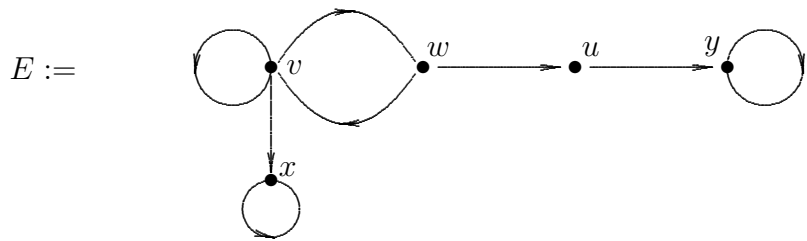
\[
\beginpicture
\setcoordinatesystem units <1cm,0.75cm>

\setplotarea x from 0 to 7, y from -3.5 to 1.6




\put{$E :=$} at 0 0

\put{$\bullet$} at 3 0

\put{$\bullet$} at 5 0

\put{$\bullet$} at 3 -2

\put{$\bullet$} at 7  0

\put{$\bullet$} at 9 0






\put{$v$}[l] at 3.2 0

\put{$w$}[b] at 5.2 0.2

\put{$x$}[l] at 3.1 -1.8

\put{$u$}[b] at 7.2 0.2

\put{$y$}[b] at 8.8 0.2
\setquadratic

\plot 3.1 0.1  4 1 4.9 0.1 /

\plot 3.1 -0.1 4 -1  4.9 -0.1 /

\circulararc 360 degrees from 3 0 center at 2.5 0

\circulararc 360 degrees from 3 -2 center at 3 -2.5

\circulararc 360 degrees from 9 0 center at 9.5 0

\arrow <0.25cm> [0.1,0.3] from 2.015 0.1 to 2 -0.1 

\arrow <0.25cm> [0.1,0.3] from 4.1 -0.985 to 3.9 -1

\arrow <0.25cm> [0.1,0.3] from 3.9 0.985 to 4.1 1

\arrow <0.25cm> [0.1,0.3] from 9.985 0.1 to 10 -0.1

\arrow <0.25cm> [0.1,0.3] from 3 -0.15 to 3 -1.85

\arrow <0.25cm> [0.1,0.3] from 2.95 -2.98 to 3.04 -3.01

\arrow <0.25cm> [0.1,0.3] from 5.2 0 to 6.8 0

\arrow <0.25cm> [0.1,0.3] from 7.2 0 to 8.8 0


\endpicture
\]

\noindent
whose proper communication sets are  $PC_1 = \{v,w\}$, $PC_2 = \{x\}$ and $PC_3 = \{y\}$.
The proper communication graph $PC(E)$ is as shown below
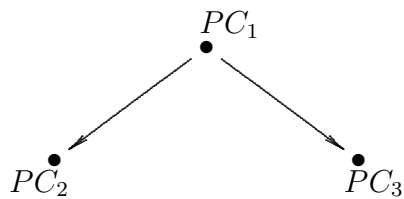
\[
\beginpicture
\setcoordinatesystem units <1cm,0.75cm>

\setplotarea x from -3 to 3, y from -3 to 1




\put{$\bullet$} at 0 0

\put{$\bullet$} at -2 -2

\put{$\bullet$} at 2 -2

\put{$PC_1$} at 0.3 0.4

\put{$PC_2$} at -2.2 -2.4

\put{$PC_3$} at 2.2 -2.4

\arrow <0.25cm> [0.1,0.3] from -0.2 -0.2 to -1.8 -1.8

\arrow <0.25cm> [0.1,0.3] from 0.2 -0.2 to 1.8 -1.8

\endpicture
\]
\end{enumerate}
\end{examples}

\noindent We now turn our attention to the effect of symbol expansion on the left Krieger cover
of a shift space. Loosely speaking, the effect of a symbol expansion $\mathcal{E}$ on a shift space $\textsf{X}$ is to 
form a new shift space $\textsf{X}'$ by adding a new letter (which we refer to as $\ast$) 
after all occurrences of certain symbol (in this case $a$) in $\textsf{X}$. 

More specifically, let $\textsf{X}$ be a shift space over the alphabet $\mathcal{A}$ and fix $a \in \mathcal{A}$. 
Choose a symbol $\ast \not\in \mathcal{A}$ and set $\mathcal{A}' = \mathcal{A} \cup \{ \ast \}$. Let
$\mathcal{E}$ be the symbol expansion map $\mathcal{E} : \mathcal{B} ( \textsf{X} ) \to ( \mathcal{A}' )^*$
which adds the symbol $\ast$ after every instance of $a$ in a word $w \in \mathcal{B} ( \textsf{X} )$.
The language $\mathcal{E} ( \mathcal{B} ( \textsf{X} ) )$ is extendible, however it is not factorisable as it does
not contain any word beginning with $\ast$. Let $\mathcal{L}$ 
be the smallest submonoid of $( \mathcal{A}' )^*$ which contains all subwords of the elements of
$\mathcal{E} ( \mathcal{B} ( \textsf{X} ) )$, then $\mathcal{L}$ is an extendible, factorisable 
language and so uniquely defines a shift space $\textsf{X}'$ by \cite[Proposition 1.3.4]{lm}.
By abuse of notation we shall also use $\mathcal{E}$ to refer to the induced maps from 
$\textsf{X}$ to $\textsf{X}'$ and $\textsf{X}^+$ to $(\textsf{X}')^+$.

Our first result shows how the symbol expansion $\mathcal{E}$ affects the connectivity of
the left Krieger covers of $\textsf{X}$ and $\textsf{X}'$.

\begin{lem} \label{lem:preservepropercommunication}
Let $\textsf{X}$, $\mathcal{E}$, $\textsf{X}\,'$ be as above.
Then for each $x^+,y^+ \in \textsf{X}^+$ there is a path labelled $w$ of length $|w| \geq 1$
from $P_\infty(x^+)$ to $P_\infty(y^+)$ in the left Krieger cover of $\textsf{X}$ if and only if there is a path
labelled $\mathcal{E}(w)$ of length $|\mathcal{E}(w)| \geq 1$ from $P_\infty(\mathcal{E}(x^+))$ to $P_\infty(\mathcal{E}(y^+))$
in the left Krieger cover of $\textsf{X}\,'$.
\end{lem}

\begin{proof} Suppose there is a path labelled $w$ of length $|w| \geq 1$ from $P_\infty(x^+)$ to $P_\infty(y^+)$ in the left Krieger cover
of $\textsf{X}$.  Then there is a right infinite ray $z^+ \in P_\infty(y^+)$ such that $x^+ = wz^+$.  But then
$\mathcal{E}(x^+) = \mathcal{E}(w)\mathcal{E}(z^+)$ and so there is a path labelled $\mathcal{E}(w)$
from $P_\infty(\mathcal{E}(x^+))$ to $P_\infty(\mathcal{E}(z^+))$ in the left Krieger cover of $\textsf{X}'$.
 It is straightforward  to check that $P_\infty(\mathcal{E}(z^+)) = P_\infty(\mathcal{E}(y^+))$ in the left
Krieger cover of $\textsf{X}'$.  Moreover since $|\mathcal{E}(w)| \geq |w|$ we certainly have $|\mathcal{E}(w)| \geq 1$.

Conversely, suppose there is a path labelled $v$ of length $|v| \geq 1$ from $P_\infty(\mathcal{E}(x^+))$ to $P_\infty(\mathcal{E}(y^+))$ in the left Krieger cover of $\textsf{X}'$.  Then there is a $z^+ \in P_\infty(\mathcal{E}(y^+))$ with $\mathcal{E}(v)z^+ = x^+$.  Note that $z^+$ cannot begin with $\ast$
since all infinite pasts of such a right-infinite ray must end in $a$.  Thus $z^+ = \mathcal{E}(t^+)$
for some $t^+ \in \textsf{X}^+$.  Now $v$ cannot begin with $\ast$ since $\mathcal{E}(x^+)$ cannot begin
with $\ast$.  Also, $v$ cannot end in $a$ since $\mathcal{E}(t^+)$ cannot begin with $\ast$.  The construction
of the left Krieger cover of $\textsf{X}'$ ensures that every $a$ which appears in $v$ is immediately followed
by $\ast$.  Thus $v = \mathcal{E}(w)$ for some $w \in {\mathcal B}(\textsf{X})$.  Since $v$ cannot
  begin in $\ast$ we must have $|w| \geq 1$.  It is then straightforward
to check that $x^+ = w t^+$ and it follows, by construction of the left Krieger cover of $\textsf{X}$,
since $P_\infty(t^+) = P_\infty(y^+)$ that there is a path labelled $w$ from $P_\infty(x^+)$ to
$P_\infty(y^+)$, establishing our result.
\end{proof}

\noindent We now give the relationship between the proper communication
graphs of the left Krieger covers of $\textsf{X}$ and $\textsf{X}\,'$.

\begin{prop}\label{prop:symexpand}
Let $\textsf{X}$, ${\mathcal E}$, $\textsf{X}\,'$ be as above.
Then there is a one-to-one correspondence between the proper communication sets
of the left Krieger covers $(E_K(\textsf{X}),{\mathcal L}_K)$ and 
$(E_K(\textsf{X}\,'),{\mathcal L}_K)$ of $\textsf{X}$ and $\textsf{X}\,'$ respectively.  Moreover, the
proper communication graphs $PC(E_K(\textsf{X}))$ and $PC(E_K(\textsf{X}\,'))$
are isomorphic as directed graphs.
\end{prop}

\begin{proof}
We note first that
\[
\left(\textsf{X}'\right)^+=\left({\mathcal E}(\textsf{X})\right)^+=\mathcal{E}(\textsf{X}^+)
\cup \sigma\left(\mathcal{E}(a\textsf{X}^+\cap \textsf{X}^+)\right)
\]
since $\mathcal{E}(\textsf{X}^+)$ consists precisely of
those one-sided infinite sequences in $\textsf{X}'$ which do not begin with
$\ast$.  By Lemma \ref{lem:preservepropercommunication} vertices  $P_\infty(x^+)$ and $P_\infty(y^+)$
in the left Krieger cover $(E_K(\textsf{X}),{\mathcal L}_K)$ of $\textsf{X}$ properly communicate if and only if the vertices $P_\infty(\mathcal{E}(x^+))$
and $P_\infty(\mathcal{E}(y^+))$ properly communicate in the left Krieger cover $(E_K(\textsf{X}'),{\mathcal L}_K)$ of $\textsf{X}'$.

We now show that vertices of the form $P_\infty(\sigma(\mathcal{E}(ax^+)))$ in $(E_K(\textsf{X}'),{\mathcal L}_K)$ do not give rise to any additional proper communication sets. In $(E_K(\textsf{X}'),{\mathcal L}_K)$ these vertices only appear in the configuration
\[
\beginpicture
\setcoordinatesystem units <2cm,0.75cm>

\setplotarea x from 0 to 6, y from -0.2 to 1.6




\put{$P_\infty(\mathcal{E}(ax^+))$}[l] at 0.3 0

\put{$P_\infty(\sigma(\mathcal{E}(ax^+)))$} at 3.2 0

\put{$P_\infty(\mathcal{E}(x^+))$} at 5.4 0

\put{$\ast$}[b] at 4.4 0.2

\put{$a$}[b] at 2 0.2



\arrow <0.25cm> [0.1,0.3] from 4 0 to 4.8 0 

\arrow <0.25cm> [0.1,0.3] from 1.5 0 to 2.4 0


\endpicture
\]
and since every past of $\sigma(\mathcal{E}(ax^+))$ necessarily ends in an
$a$, we see that these two edges are the only ones emitted and
received by $P_\infty(\sigma(\mathcal{E}(ax^+)))$.

This implies that if  $P_\infty(x^+)$ and $P_\infty(ax^+)$ properly communicate
in the $(E_K(\textsf{X}),{\mathcal L}_K)$, then $P_\infty(\mathcal{E}(x^+))$, $P_\infty(\mathcal{E}(ax^+))$,
and $P_\infty(\sigma(\mathcal{E}(ax^+))$ properly communicate in $(E_K(\textsf{X}'),{\mathcal L}_K)$ and so $P_\infty(\mathcal{E}(x^+))$, $P_\infty(\mathcal{E}(ax^+)$ and $P_\infty(\mathcal{E}(\sigma(\mathcal{E}(ax^+)))$ are all elements of the same proper communication set.

   On the other hand, if there is no path from
$P_\infty(x^+)$ to $P_\infty(ax^+)$ in $(E_K(\textsf{X}),{\mathcal L}_K)$ there can be no path from
$P_\infty(\mathcal{E}(x^+))$ to $P_\infty(\sigma(\mathcal{E}(ax^+)))$ and so these elements do not properly communicate.   Consequently,  symbol expansion does not change the number of proper communication sets.

Let $PC_i^\textsf{X}$ denote the proper communication sets of $(E_K(\textsf{X}),{\mathcal L}_K)$ and
let $PC_i^{\textsf{X}'}$ denote the proper communication sets of $(E_K(\textsf{X}'),{\mathcal L}_K)$.  Let $\phi_{PC}^0$ denote the bijection between the proper communication sets of $(E_K(\textsf{X}),{\mathcal L}_K)$ and $(E_K(\textsf{X}'),{\mathcal L}_K)$ which satisfies the condition that for all properly communicating vertices $P_\infty(x^+)$ in $(E_K(\textsf{X}),{\mathcal L}_K)$ we have $P_\infty(x^+) \in PC^\textsf{X}_i$ if and only if  $P_\infty(\mathcal{E}(x^+)) \in \phi_{PC}^0(PC^{\textsf{X}}_i)$.  It follows from Lemma \ref{lem:preservepropercommunication} that there is an
edge with $s(e) = PC^{\textsf{X}}_i$ and $r(e) = PC^{\textsf{X}}_j$ in $PC(E_K(\textsf{X}))$ if and only if there is an edge $\phi_{PC}^1(e)$ with $s(\phi_{PC}^1(e)) =
\phi_{PC}^0(PC^{\textsf{X}}_i)$ and $r(\phi_{PC}^1(e)) = \phi_{PC}^0(PC^{\textsf{X}}_j)$ in $PC(E_K(\textsf{X}'))$.  Thus the map
$\phi_{PC} = (\phi^0_{PC},\phi^1_{PC}) : PC(E_K(\textsf{X})) \to PC(E_K(\textsf{X}'))$ is a directed graph
isomorphism and our result is established.
\end{proof}

\noindent For sofic shifts we can go even further:

\begin{thm}\label{thm:flowinv}
The proper communication graph $PC(E_K(\textsf{X}))$ is a flow invariant for sofic shifts.
\end{thm}

\begin{proof}
By Proposition \ref{prop:symexpand}  the proper communication graph is preserved under symbol expansion.
It is shown in \cite{kri} that the left Krieger cover is a conjugacy invariant for sofic shifts, and so
the proper communication graph will also be preserved under conjugacy.  Our result follows since
symbol expansion and conjugacy are the operations which generate flow equivalence for shifts by \cite{parsu}.
\end{proof}

\begin{cor}
Let $\textsf{X}$ be a sofic shift whose left Krieger cover is reducible. Suppose that $\textsf{X}'$ is flow equivalent to
$\textsf{X}$, then $\textsf{X}'$ is sofic with a reducible left Krieger cover.
\end{cor}

\begin{proof}
By \cite[Corollary 3.2.3]{lm} any shift conjugate to a sofic shift is also sofic. Suppose that $\textsf{X}'$ is obtained
from $\textsf{X}$ by a symbol expansion $\mathcal{E}$. The proof of Proposition \ref{prop:symexpand} shows that
the left Krieger cover of $\textsf{X}'$ is obtained from that of $\textsf{X}$ by adding a vertex $P_\infty ( \sigma ( \mathcal{E} ( a x^+ ) ) )$
between $P_\infty ( \mathcal{E} ( ax^+ ) )$ and $P_\infty ( \mathcal{E} ( x^+ ) )$ with incoming label $a$ and outgoing
label $\ast$ for every vertex of the form $P_\infty ( a x^+ )$ in the left Krieger cover of $\textsf{X}$. The resulting labelled graph will be finite and so the shift $\textsf{X}'$ will be sofic, and the result follows from \cite{parsu} once again.

The remaining statement follows immediately from Theorem \ref{thm:flowinv}.
\end{proof}

\section{Reducibility of the left Krieger cover of an AFT shift} \label{sec:aftresult}

\noindent The following example is adapted from the one found in \cite[\S 2.3 Figure 2]{bfp}.

\begin{example} \label{ex:bfg}
Consider the sofic shift $\textsf{X}$ with left Fischer cover as shown
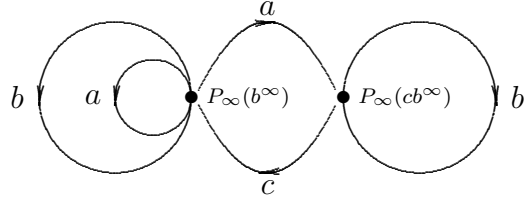
\[
\beginpicture

\setcoordinatesystem units <1cm,1cm>

\setplotarea x from 0 to 8, y from 0 to 1.5

\put{$\bullet$} at 3 0

\put{$\bullet$} at 5 0

\put{{\tiny $P_\infty(b^\infty)$}}[l] at 3.2 0

\put{{\tiny $P_\infty(cb^\infty)$}}[l] at 5.2 0

\put{$a$}[l] at 1.6 0

\put{$a$}[b] at 4 1.1

\put{$c$}[t] at 4 -1.1

\put{$b$}[r] at 0.8 0

\put{$b$}[l] at 7.2 0


\setquadratic

\plot 3.1 0.1  4 1 4.9 0.1 /

\plot 3.1 -0.1 4 -1  4.9 -0.1 /

\circulararc 360 degrees from 3 0 center at 2.5 0

\circulararc 360 degrees from 5 0 center at 6 0

\circulararc 360 degrees from 3 0 center at 2 0

\arrow <0.25cm> [0.1,0.3] from 2.015 0.1 to 2 -0.1 

\arrow <0.25cm> [0.1,0.3] from 4.1 -0.985 to 3.9 -1

\arrow <0.25cm> [0.1,0.3] from 3.9 0.985 to 4.1 1

\arrow <0.25cm> [0.1,0.3] from 1.015 0.1 to 1 -0.1

\arrow <0.25cm> [0.1,0.3] from 7 0.1 to 7.015 -0.1

\endpicture
\]

\noindent
The shift $\textsf{X}$ is irreducible by Remark \ref{rem:soficirred} and is strictly sofic since
the left Fischer cover has two distinct edges with the same label.
Since there are two distinct representatives of the right-ray $ab^\infty$
starting at the same vertex, the left Fischer cover is not right-closing which implies
that $\textsf{X}$ is not AFT by Theorem \ref{thm:aftconditions}.
One checks that any right-ray $x^+$ which only contains the symbols $a, b$ has
predecessor set  $P_\infty(x^+) = P_\infty ( b^\infty ) = \textsf{X}^-$. Any right-ray of the
form $wc y$ where $y \in \textsf{X}^+$ and $w$ contains an $a$ also has predecessor
set $P_\infty(wcy) =  \textsf{X}^- = P_\infty ( b^\infty)$. However any right-ray of the form $b^n c y$ where $y \in \textsf{X}^+$ and $n \ge 0$
is such that $z =\ldots (ac)(ac) \in \textsf{X}^-$ does not belong to $P_\infty ( b^n cy )$.
Since every representative of the right-ray $b^n cy$ must start at the same vertex we conclude that
$P_\infty ( b^n cy ) = P_\infty ( c b^\infty )$
for all $n \ge 0$ and so the left Krieger and Fischer covers coincide. In particular,
the left Krieger cover is irreducible.
\end{example}

On the other hand, Example \ref{ex:evenshift} shows that there are strictly sofic shifts
with reducible left Krieger covers. It is then natural to try to identify the class of strictly sofic shifts for which the
left Krieger cover is reducible. In this section we produce a partial answer to this question in
Theorem \ref{thm:aftresult}: If the shift is AFT then the left Krieger cover is reducible.
However, by Example \ref{ex:jonathanstrikesback} we see that there are irreducible non AFT shifts
which have a reducible left Krieger cover.

Before giving the proof of our main result we need a few technical results concerning
the left Krieger cover of an AFT shift.
By Remark \ref{rem:lfc} (ii) we shall identify the left Fischer cover of an irreducible
strictly sofic shift $\textsf{X}$ with the minimal irreducible subgraph of the left
Krieger cover of $\textsf{X}$.

The following result shows that the vertices in the left Krieger cover of an irreducible strictly
sofic shift that correspond to non intrinsically synchronising predecessor sets
 do not connect to the vertices of the left Fischer cover (see Remark \ref{rem:lfc} (ii)).

\begin{lem} \label{lem:nottrans}
Let $\textsf{X}$ be an irreducible strictly sofic shift.
Suppose $\textsf{X}$ has a predecessor set $P_\infty(x^+)$ that is not instrinsically synchronising.
Then the left Krieger cover $(E_K,{\mathcal L}_K)$ of $\textsf{X}$ is reducible.
\end{lem}

\begin{proof}
We show that there is no intrinsically synchronising predecessor set $P_\infty(y^+)$
and path $\lambda$ in $E_K^*$ with $s(\lambda) = P_\infty(x^+)$ and $r(\lambda) =
P_\infty(y^+)$.  Suppose otherwise, then without loss of generality we may assume that $y^+$ is
an intrinsically synchronising ray for $\textsf{X}$.  By definition of the
left Krieger cover we must have $P_\infty (\mathcal{L}_K ( \lambda ) y^+ ) = P_\infty (x^+)$.
Since $y^+$  is an intrinsically synchronising ray,
it follows that $P_\infty(x^+)$ is an intrinsically synchronising predecessor set,
a contradiction.  Thus the left Krieger cover of $\textsf{X}$ is reducible.
\end{proof}

\noindent The following lemma gives a characterisation of \is \ right-rays in a sofic shift in
terms of their representatives in the left Krieger cover.

\begin{lem} \label{lem:isvertsareunique}
Let $\textsf{X}$ be a sofic shift, then $w \in \mathcal{B} (\textsf{X})$ is
\is  \ if and only if $P_\infty ( w x^+ ) = P_\infty ( w y^+ )$ for all $x^+ , y^+ \in \textsf{X}^+$
such that $w x^+ , w y^+  \in \textsf{X}^+$. In particular, all representatives of $w$ in the
left Krieger cover of $\textsf{X}$ begin at the same vertex.
\end{lem}

\begin{proof}
Suppose that $w \in \mathcal{B} ( \textsf{X} )$ is an \is \ word, and that $x^+ , y^+ \in \textsf{X}^+$
are such that $P_\infty ( w x^+ ) \neq P_\infty ( w y^+ )$. Then, without loss of generality, we may
 assume that there is $x^- \in P_\infty ( w x^+ )$
such that $x^- w y^+ \not\in \textsf{X}$. Since $x^- w x^+ \in \textsf{X}$, $w y^+ \in \textsf{X}^+$
and $x^- w y^+ \not\in \textsf{X}$ it follows from \cite[Corollary 1.3.5]{lm} that there are $m , n \geq 1$ such that
$x^-_{-n} \ldots x^-_{-1} w \in \mathcal{B} ( \textsf{X} )$ and $w y_1 \ldots y_m \in \mathcal{B} ( \textsf{X} )$
but $x^-_{-n} \ldots x^-_{-1} w y_1 \ldots y_m \not\in \mathcal{B} ( \textsf{X} )$. But this contradicts
the hypothesis that $w$ is an \is \ word, hence $P_\infty ( w x^+ ) \subseteq P_\infty ( w y^+ )$. A symmetric
argument shows that $P_\infty ( w y^+ ) \subseteq P_\infty ( w x^+ )$ and hence $P_\infty ( w x^+ ) =
P_\infty  ( w y^+ )$ as required.

Conversely, suppose that $w \in \mathcal{B} ( \textsf{X} )$ is such that
$P_{\infty} ( w x^+ ) = P_\infty (w y^+ )$ for all $x^+ , y^+ \in \textsf{X}^+$
such that $w x^+ , w y^+  \in \textsf{X}^+$. Then by Remark \ref{rem:lkc} (i) every path in
$(E_K ,{\mathcal L_K})$ labelled $w$ begins
at $P_\infty(w x^+)$.  If $uw$ and $wv \in {\mathcal B} ( \textsf{X} )$ then there are paths
$\mu \nu$ and $\nu' \lambda \in E_F^*$ with $|\nu| = |\nu'| = |w|$ and
${\mathcal L_F} ( \mu \nu ) = uw$ and ${\mathcal L}_F( \nu' \lambda ) = wv$.
Since ${\mathcal L}_F( \nu ) = {\mathcal L}_F( \nu' ) = w$,  by hypothesis we must have
$s ( \nu ) = s( \nu' ) = P_\infty (w x^+) = r(\mu)$ and so $\mu \nu' \lambda \in E_F^*$.
Thus ${\mathcal L}_F( \mu \nu' \lambda ) = uwv \in {\mathcal B } ( \textsf{X} )$ and
so $w$ is an intrinsically synchronising word.
\end{proof}

\noindent
The final statement of Lemma \ref{lem:isvertsareunique} is an analogue of \cite[Lemma 1.1]{trow}
for left Krieger covers.



\noindent
For completeness we include a proof of the following Lemma (cf. \cite[Lemma 4.2]{kri}).

\begin{lem} \label{lem:onex}
Let $\textsf{X}$ be a strictly sofic shift. Then there is a word $w \in \mathcal{B} ( \textsf{X} )$
such that $x = \ldots w w \ldots \in \textsf{X}$ has more than one representative
in the left Fischer cover of $\textsf{X}$.  In particular there are at least two distinct circuits,
$\alpha, \beta \in E_F^*$ with ${\mathcal L}_F ( \alpha ) = {\mathcal L}_F( \beta ) = w$.
\end{lem}

\begin{proof}
Let $\textsf{X}$ be a strictly sofic shift over $\mathcal{A}$. Then the factor map
$\pi_{\mathcal L} : \textsf{X}_{E_F} \to \textsf{X}$
induced by the $1$-block map ${\mathcal L}_F :E^1_F \to {\mathcal A}$ fails to be injective at some point $x \in
\textsf{X}$.  Let $y,z \in \pi_{\mathcal{L}}^{-1}(x)$ with $y \ne z$.  We claim that there is an $n \in {\mathbf Z}$
such that $\{(r(y_i),r(z_i)) : i \ge n \} \subseteq E_F^0 \times E_F^0$ does not contain $(v,v)$ for any
$v \in E_F^0$.   Since $y \ne z$ there is an $n \in {\mathbf Z}$ such that $y_n \ne z_n$.  Let $u_n = r(y_n)$ and
$v_n = r(z_n)$.  Then since ${\mathcal  L}(y_n) = {\mathcal L} (z_n)$ and the left Fischer cover is left-resolving
we cannot have $u_n = v_n$.  Now the edges $y_{n + 1}$ and $z_{n+1}$ begin at different vertices so $y_{n+1} \ne z_{n+1}$.
An inductive argument completes the proof of the claim.

As $E^0_F \times E^0_F$ is finite there must be $(u,v) \in E^0_F \times E^0_F$ with $u \ne v$ and $m \ge 1$
such that $(r(y_i),r(z_i)) = (u,v) = (r(y_{i + m}), r(z_{i+m}))$ where $i \ge n$.  Hence
$y_{i+1} \dots y_{i+m}$ and $z_{i+1} \dots z_{i+m}$ are distinct circuits $\alpha$ and $\beta$ in $E_F$
with ${\mathcal L}_F(\alpha) = {\mathcal L}_F(\beta) = w$, say.  Hence $\dots \alpha \alpha \dots $
and $\dots \beta \beta \dots$ are the required representatives of $x$ in the left Fischer cover.
\end{proof}

\begin{lem} \label{lem:twox}
Let $\textsf{X}$ be a strictly sofic shift and $w \in \mathcal{B} ( \textsf{X} )$ be as in Lemma \ref{lem:onex}.
Then $w$ is not an intrinsically synchronising word.
Moreover, for every positive integer $k_0$ there is an integer $k$ with
$k > k_0$ and $x^+ \in \textsf{X}^+$ such that $P_\infty ( w^{k_0}x^+ ) =P_\infty ( w^k x^+ )$.
\end{lem}

\begin{proof}
Now by Lemma \ref{lem:onex} there are circuits
$\alpha \neq \beta$ in $E_F^*$ with $r(\alpha) = s(\alpha) \neq r(\beta) = s(\beta)$ such that
${\mathcal L}_F(\alpha) = {\mathcal L}_F(\beta) = w$.  If $w$ is intrinsically synchronising,
then Lemma \ref{lem:isvertsareunique} implies that
$s(\alpha) = s(\beta) = P_\infty (w x^+)$ for some $x^+ \in \textsf{X}^+$, a contradiction.

Since $\textsf{X}$ is a sofic shift there are a finite number of distinct predecessor sets.
It follows, by the pigeonhole principle, that for each positive integer $k_0$ we can find an
integer $k$ with $k  > k_0$ and an $x \in \textsf{X}^+$ such that
$P_\infty (w^{k_0}x^+) = P_\infty (w^k x^+)$ as required.
\end{proof}

\noindent
We now prove that the left Krieger cover of a strictly sofic, irreducible AFT shift $\textsf{X}$ is reducible.

\begin{thm}  \label{thm:aftresult}
Let $\textsf{X}$ be an irreducible, strictly sofic AFT shift.  Then the left Krieger cover $(E_K,{\mathcal L}_K)$
of $\textsf{X}$ is reducible.
\end{thm}

\begin{proof}
Suppose, for contradiction that the left Fischer cover and the left Krieger cover of $\textsf{X}$ coincide,
in particular $E_K = E_F$ is irreducible by Remarks \ref{rem:lfc} (ii).
Let $w$ be as in Lemma \ref{lem:onex}, let $k_0$ be a positive integer and let $k > 0$ and $x^+ \in \textsf{X\,}^+$
satisfy the conditions of Lemma \ref{lem:twox} with respect to $w$ and $k_0$.
Let $\alpha$ and $\beta$ be distinct circuits in $E_F^*$ representing  $w$ as in Lemma \ref{lem:onex}.
It follows by Remark \ref{rem:lkc} (i) there is a path $\mu \in E_F^*$ from $P_\infty (w^k x^+)$ to
$P_\infty (w^{k_0}x^+)$ labelled $w^{k - k_0}$.
Since $P_\infty(w^k x^+) = P_\infty (w^{k_0} x^+)$, $\mu$ is a circuit.  Note that
since $\alpha$ and $\beta$ are distinct representatives of $w$, $\mu$ must be distinct from one of $\alpha^{k-k_0}$
and $\beta^{k-k_0}$.  Without loss of generality we may assume that $\mu$ is distinct from $\beta^{k-k_0}$.

Since the left Krieger cover and the left Fischer cover coincide the predecessor set $P_\infty (w^{k_0} x^+)$ is \is, say
$P_\infty ( w^{k_0} x^+ ) = P_\infty (m y^+)$ where
$m \in \mathcal{B} ( \textsf{X} )$ is an \is \ word.  By Remark \ref{rem:lkc} (i) there is a path $\nu$
starting at $P_\infty( m y^+)$ with ${\mathcal L}_F ( \nu ) = m$.
Since $E_F$ is irreducible there is a path $\lambda$ from $r(\nu)$ to $s(\beta)$.
Let $v = {\mathcal L}_F(\lambda)$.  Then $mvw$ is represented by
$\nu\lambda \beta$.  Since $m$ is intrinsically synchronising we have
$P_\infty (mvw^{k_0} x^+) = P_\infty (m y^+)$ by Lemma \ref{lem:isvertsareunique}.
By Remark \ref{rem:lkc} (i) there is a path $\kappa \in E_F^*$ from $P_\infty ( mvw^{k_0} x^+)$
to $P_\infty (w^{k_0} x^+)= P_\infty(m y^+)$ labelled by $mv$.
Since $P_\infty (mvw^{k_0} x^+) = P_\infty (m y^+)$ it follows that $\kappa$
is a circuit. Since $(E_F,\mathcal{L}_F)$ is left-resolving and
$\mathcal{L}_F ( \kappa ) = mv = \mathcal{L}_F ( \nu\lambda)$ we must have
$r ( \kappa ) \neq r ( \nu \lambda )$ and hence $\kappa \ne \nu\lambda $.

Since $\textsf{X}$ is AFT, $\pi_{\mathcal L}$ is right-closing with delay $D$ by
Theorem \ref{thm:aftconditions} (2).  By choosing $n$ large enough so that
$| mvw^{n(k - k_0)} | \ge D + | m v | + 1$ we obtain $\kappa = \nu \lambda$
by Remark \ref{rem:delaypath}, a contradiction.  Hence our assumption that the left
Krieger cover was equal to the left Fischer cover must have been false. In particular,
the predecessor set $P_\infty ( w^{k_0} x^+ )$ is
not \is\ and the result then follows from Lemma \ref{lem:nottrans}.
\end{proof}

\noindent We are now able to deduce the following result concerning the
Matsumoto algebra associated to a shift space (for more details see \cite{cm,m}).
Our argument employs graph $C^*$-algebra techniques.
For more details about graph $C^*$-algebras see \cite{bprsz}.

\begin{cor} \label{cor:notsomatso}
Let $\textsf{X}$ be an irreducible strictly sofic AFT shift, then the
Matsumoto algebra $\mathcal{O}_{\textsf{X}}$ associated to $\textsf{X}$
is not simple.
\end{cor}

\begin{proof}
By \cite[Corollary 6.8]{bp2} (see also \cite[Theorem 3.5]{c}, \cite[Proposition 7.1]{m5}) $\mathcal{O}_{\textsf{X}}$ is
isomorphic to the graph algebra  $C^* ( E_K )$.
By Theorem \ref{thm:aftresult} $E_K$ is reducible and
the result then follows by \cite[Proposition 5.1]{bprsz}.
\end{proof}

\begin{rmk} \label{rmk:infinitepastfinitefuture}  Another left-resolving cover of a sofic shift $\textsf{X}$ over the
alphabet $\mathcal{A}$ may be defined as follows.
For $w \in {\mathcal B}(\textsf{X})$ let
$P_\infty^f ( w ) = \{ x^- \in \textsf{X}^- \,:\, x^- w \in \textsf{X}^-\}$ be the collection
of all left-rays which may precede $w$.  The cover is given by the labelled graph
$( E_\infty^f , {\mathcal L} )$ where the vertices of $E_\infty^f$ are the predecessor
sets $P_\infty^f ( w)$ and there is an edge labelled
$a \in {\mathcal A}$ from $P_\infty^f ( w )$ to $P_\infty^f ( v )$ if and only if  $a v \in {\mathcal B}( \textsf{X})$ and
$P_\infty^f ( w) = P_\infty^f ( a v )$.  One may check that the labelling ${\mathcal L}$ is well-defined.

A corresponding version of Theorem \ref{thm:aftresult} may be proved: Let $\textsf{X}$ be an irreducible, strictly
sofic AFT shift.  Then the cover $( E_\infty^f , {\mathcal L} )$ of $\textsf{X}$ is reducible.
Since the sets $P_\infty^f ( w)$ consist of left-infinite rays, the argument used to prove Theorem
\ref{thm:aftresult} applies \emph{mutatis mutandis}.
\end{rmk}

\noindent
The converse of Theorem \ref{thm:aftresult} does not hold as we see in the following example which was inspired
by \cite[Figure 7.5.1]{sam}.

\begin{example} \label{ex:jonathanstrikesback}
Consider the sofic shift $\textsf{Y}$ with left Fischer cover as shown below
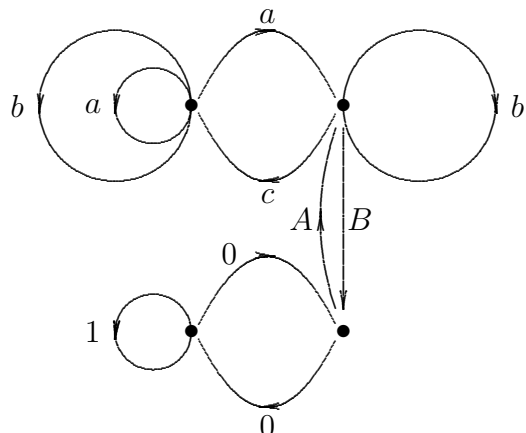
\[
\beginpicture

\setcoordinatesystem units <1cm,1cm>

\setplotarea x from 0 to 8, y from -4 to 1.5

\put{$\bullet$} at 3 0

\put{$\bullet$} at 5 0



\put{$a$}[l] at 1.6 0

\put{$a$}[b] at 4 1.1

\put{$c$}[t] at 4 -1.1

\put{$b$}[r] at 0.8 0

\put{$b$}[l] at 7.2 0


\setquadratic

\plot 3.1 0.1  4 1 4.9 0.1 /

\plot 3.1 -0.1 4 -1  4.9 -0.1 /

\circulararc 360 degrees from 3 0 center at 2.5 0

\circulararc 360 degrees from 5 0 center at 6 0

\circulararc 360 degrees from 3 0 center at 2 0

\arrow <0.25cm> [0.1,0.3] from 2.015 0.1 to 2 -0.1 

\arrow <0.25cm> [0.1,0.3] from 4.1 -0.985 to 3.9 -1

\arrow <0.25cm> [0.1,0.3] from 3.9 0.985 to 4.1 1

\arrow <0.25cm> [0.1,0.3] from 1.015 0.1 to 1 -0.1

\arrow <0.25cm> [0.1,0.3] from 7 0.1 to 7.015 -0.1

\put{$\bullet$} at 3 -3

\put{$\bullet$} at 5 -3



\put{$1$}[l] at 1.6 -3

\put{$0$}[b] at 3.5 -2.1

\put{$0$}[t] at 4 -4.1


\setquadratic

\plot 3.1 -2.9  4 -2 4.9 -2.9 /

\plot 3.1 -3.1 4 -4  4.9 -3.1 /

\plot 4.9 -0.3  4.7 -1.5  4.9 -2.7 /

\circulararc 360 degrees from 3 -3 center at 2.5 -3

\arrow <0.25cm> [0.1,0.3] from 2.015 -2.9 to 2 -3.1 

\arrow <0.25cm> [0.1,0.3] from 3.9 -1.985 to 4.1 -2

\arrow <0.25cm> [0.1,0.3] from 4.1 -3.985 to 3.9 -4

\arrow <0.25cm> [0.1,0.3] from 4.7 -1.5 to 4.7 -1.4


\put{$A$}[r] at 4.6 -1.5

\put{$B$}[l] at 5.05 -1.5

\arrow <0.25cm> [0.1,0.3] from 2.015 -2.9 to 2 -3.1

\arrow <0.25cm> [0.1,0.3] from 5 -0.3 to 5 -2.7

\endpicture
\]
\noindent The shift $\textsf{Y}$ is irreducible by Remark \ref{rem:soficirred}, and since the
left Fischer cover has two distinct edges with the same label it follows that $\textsf{Y}$ is strictly sofic.
Since there are two distinct representatives of the right-ray $ab^\infty$
starting at the same vertex, the left Fischer cover is not right-closing which implies
that $\textsf{X}$ is not AFT by Theorem \ref{thm:aftconditions}. However $\textsf{Y}$ contains
the even shift as a subshift.  As in Example \ref{ex:evenshift} there are three different
predecessor sets associated to right rays in the symbols $0$ and $1$.  Hence the left Krieger cover
of $\textsf{Y}$  contains a subgraph isomorphic to the left Krieger cover of the even shift and so is
reducible (cf.\ Example \ref{ex:evenshift}).
\end{example}

\section{One-sided shifts and the past set cover}\label{sec:aftresult2}

Recall from Definition \ref{def:rightrays}, that for a shift space $\textsf{X}$, the set $\textsf{X}^+$ of right-rays consists of
all sequences $x_0 x_1 \dots$ where $x \in \textsf{X}$. This space is still invariant under the shift map $\sigma$, and
the resulting pair $(\textsf{X}^+,\sigma)$ is called a one-sided shift space (cf. \cite[p140]{lm}). Where possible we
shall indicate that we are working with a one-sided shift by adding the ${}^+$ superscript. Most of  the
concepts and definitions we have used for two-sided shifts still apply to one-sided shifts.  For example, any presentation of a two-sided
shift is a presentation of the corresponding one-sided shift (simply consider right-infinite paths).
Hence the one-sided sofic shifts are the ones presented by finite labelled graphs (cf.\ \cite[Definition 3.1.3]{lm}).

For one-sided shifts the left Krieger cover construction given in Definition \ref{infinitekrieger}
does not make sense as there are no left-infinite rays.
In this situation it is usual to work with the past-set cover (see Definition \ref{def:psc}). We give the definition for two-sided shifts,
 and note that it makes sense for one-sided shifts.

\begin{dfn}
Let $\textsf{X}$ be a shift space and $w$ a word in $\mathcal{B} ( \textsf{X} )$. The
{\em predecessor set} $P_f (w)$ of $w$ in $\textsf{X}$ is the set of
all words that can precede $w$ in $\textsf{X}$; that is,
\[
P_f (w) = \{ v \in \mathcal{B} ( \textsf{X} ) : v w  \in \mathcal{B}
( \textsf{X} ) \} .
\]
\end{dfn}

\noindent By suitably adapting \cite[Theorem 3.2.10]{lm} one may show that $\textsf{X}^+$ is a one-sided sofic
shift if and only if the set of all predecessor sets $P_f(w)$  is finite.

\begin{dfn} \label{def:psc}
Suppose that $\textsf{X}$ is a shift space over $\mathcal{A}$.
The \emph{past set cover of $\textsf{X}$} is the presentation $(E_P ,{\mathcal L}_P)$
where the vertices of $E_P$ are the predecessor
sets  $P_f(w)$. For predecessor sets $P_f ( w ) ,P_f ( v )$
and $a \in \mathcal{A}$ there is an edge labelled
$a$ from $P_f ( w )$ to $P_f ( v )$ if and only if
$P_f ( a v ) = P_f ( w )$.
\end{dfn}

\noindent As in \cite[p.73]{lm} one checks that the labelling, $\mathcal{L}_P$ is well-defined
and left-resolving.

One may again define the left Fischer cover of a one-sided sofic shift $\textsf{X}^+$ to be the minimal
left-resolving cover of $\textsf{X}^+$.

\begin{lem} \label{lem:zero}
Let $\textsf{X}$ be a (one- or two-sided) sofic shift.
Let $w \in \mathcal{B} ( \textsf{X} )$. Then $w$ is intrinsically synchronising if and only
if whenever $v \in \mathcal{B} ( \textsf{X} )$ is such that $wv \in \mathcal{B} ( \textsf{X} )$
we have $P_f ( wv ) = P_f (w)$. In particular every path in $( E_P , \mathcal{L}_P )$ labelled $w$
begins at $P_f (w)$.
\end{lem}

\begin{proof}
Let $w \in \mathcal{B}( \textsf{X} )$ be such that for all
$v \in \mathcal{B}(\textsf{X} )$ with $wv \in \mathcal{B}(\textsf{X} )$ we have
$P_f (wv) = P_f (w)$.  Suppose that $u \in \mathcal{B} ( \textsf{X} )$ is such that
$uw \in \mathcal{B}(\textsf{X} )$.
Then $u \in P_f (w) = P_f (wv)$, by assumption and so we must have $uwv \in \mathcal{B} ( \textsf{X} )$.
Thus $w$ is intrinsically synchronising.

Conversely, let $m \in \mathcal{B} ( \textsf{X} )$ be intrinsically synchronising and
$v \in \mathcal{B} ( \textsf{X} )$ be such that $mv \in \mathcal{B} ( \textsf{X} )$.
Suppose that $u \in \mathcal{B} ( \textsf{X} )$ is such that $u \in P_f (m)$ then
$umv \in \mathcal{B} ( \textsf{X} )$ and so $u \in P_f (mv)$, which implies that
$P_f (m) \subseteq P_f ( mv )$. However, by definition of $P_f (mv)$ we have $P_f ( mv ) \subseteq P_f ( m )$,
and so $P_f ( mv ) = P_f ( m )$
as required.
\end{proof}

\begin{lem} \label{lem:blmagogo}
Let $\textsf{X}$ be a (one- or two-sided) irreducible sofic shift with past set cover $( E_P , \mathcal{L}_P )$.
Let $E_P^m$ denote the subgraph of $E_P$ with vertices $P_f ( m )$ where $m \in \mathcal{B} ( \textsf{X} )$
is \is \ and the edges in $E_P$ connecting them. Then $( E_P^m , \mathcal{L}_P )$ is a minimal
left-resolving irreducible presentation of $\textsf{X}$ and so is
isomorphic (as a labelled graph) to the left Fischer cover $( E_F , \mathcal{L}_F )$ of $\textsf{X}$.
\end{lem}

\begin{proof}
Let $m_1 , m_2 \in \mathcal{B} ( \textsf{X} )$ be \is , then since $\textsf{X}$ is
irreducible there is $u \in \mathcal{B} ( \textsf{X} )$ such that $m_1 u m_2 \in \mathcal{B} ( \textsf{X} )$.
Since all subwords of $m_1 u m_2$ which contain $m_2$ are also \is \ it follows that there is a path
in $E_P^m$ labelled $m_1u$ from $P_f (m_1 u m_2 ) = P_f (m_1 )$ to $P_f ( m_2 )$. Hence $E_P^m$ is irreducible.

The labelled graph $( E_P^m , \mathcal{L}_P )$ is left-resolving as it is a subgraph of $( E_P , \mathcal{L}_P )$.
We claim that $( E_P^m , \mathcal{L}_P )$ is a presentation of $\textsf{X}$. Given $u \in \mathcal{B} ( \textsf{X} )$ the irreducibility
of $\textsf{X}$ tells us that there are $m_1 , v_1 \in \mathcal{B} ( \textsf{X} )$ with
$m_1 v_1 u \in \mathcal{B} ( \textsf{X} )$ and $m_1$ \is . Irreducibility again tells us that there are
$m_2 , v_2 \in \mathcal{B} ( \textsf{X} )$ with $(m_1 v_1 u ) v_2 m_2 \in \mathcal{B} ( \textsf{X} )$ where
$m_2$ is \is. Hence there is a path labelled
$m_1 v_1 u v_2$ starting at $P_f (m_1)$ in $( E_P^m , \mathcal{L}_P )$, and $u$ is a finite labelled path
in $(E_P^m,\mathcal{L}_P)$, which establishes our claim.
If we omit the vertex $P_f (m)$ from $E_F^m$ then Lemma \ref{lem:zero} tells us that we cannot represent $m$ in
the resulting labelled graph, and so the presentation $( E_P^m , \mathcal{L}_P )$ is minimal.
The final statement follows from the left-resolving analogue of \cite[Theorem 3.3.18]{lm}.
\end{proof}

\noindent We shall henceforth identify the left Fischer cover of an irreducible sofic shift with a subgraph
of the past set cover (cf.\ Remarks \ref{rem:lfc}).
The proof of the following Lemma is a minor modification of that given for
Lemma \ref{lem:nottrans}.

\begin{lem} \label{lem:nottransx}
Let $\textsf{X}$ be an irreducible (one- or two-sided) strictly sofic shift.
Let $w \in  \mathcal{B}(\textsf{X})$ be such that $P_f( w ) \neq P_f( m )$
for any instrinsically synchronising word $m \in \mathcal{B}(\textsf{X})$.
Then the past set cover $(E_P,{\mathcal L}_P)$ of $\textsf{X}$ is reducible.
\end{lem}

The following results for one-sided shifts can be proved in a similar
way to Lemma \ref{lem:onex} and Lemma \ref{lem:twox}.

\begin{lem} \label{lem:one}
Let $\textsf{X\,}^+$ be a one-sided strictly sofic shift. Then there is $w \in \mathcal{B} ( \textsf{X\,}^+)$
such that $x = w w \ldots \in \textsf{X\,}^+$ has more than one representative
in the left Fischer cover of $\textsf{X\,}^+$.  More specifically, there are at least two distinct circuits,
$\alpha, \beta \in E_F^*$ with ${\mathcal L}_F ( \alpha ) = {\mathcal L}_F ( \beta ) = w$.
\end{lem}

\begin{lem} \label{lem:two}
Let $\textsf{X\,}^+$ be a one-sided strictly sofic shift and let $w$ be as in Lemma \ref{lem:one}.
Then $w$ is not intrinsically synchronising.
Moreover, for every positive integer $k_0$ there is a positive integer $k$ with $k > k_0$ such that
$P_f ( w^{k_0} ) = P_f ( w^k )$.
\end{lem}

\noindent
The proof of the following theorem is similar to the one we have given for Theorem \ref{thm:aftresult}.

\begin{thm} \label{thm:aftresult2}
Let $\textsf{X}$ be a (one- or two-sided) irreducible, strictly sofic AFT shift.
Then the past set cover $( E_P ,{\mathcal L}_P )$ of $\textsf{X}$ is reducible.
\end{thm}

\begin{rmks}\label{rmk:othercover}
\begin{enumerate}
\item  Another left-resolving cover of a (one- or two-sided) sofic shift $\textsf{X}$
over the alphabet $\mathcal{A}$ may be defined as follows.
For $x^+ \in \textsf{X}^+$ let
\[
P^\infty_f ( x^+ ) = \{ w \in {\mathcal B}(\textsf{X}) \;:\; w x^+ \in \textsf{X}^+\}
\]

\noindent be the collection
of words which may precede $x^+$.  The cover is given by the labelled graph
$( E^\infty_f , {\mathcal L} )$ where the vertices of $E^\infty_f$ are the predecessor
sets $P^\infty_f ( x^+ )$ and there is an edge labelled
$a \in \mathcal{A}$ from $P^\infty_f ( x^+ )$ to $P^\infty_f ( y^+ )$ if and only if  $a y^+ \in  \textsf{X}^+$ and
$P^\infty_f ( x^+ ) = P^\infty_f ( a y^+ )$.  One checks that the labelling ${\mathcal L}$ is well-defined.

A corresponding version of Theorem \ref{thm:aftresult2} may be proved: Let $\textsf{X}$ be an irreducible,
(one- or two-sided) strictly sofic AFT shift.  Then the cover $( E^\infty_f , {\mathcal L} )$ of $\textsf{X}$ is reducible.
Since the sets $P^\infty_f ( x^+ )$ consist of elements of ${\mathcal B}(\textsf{X})$, the argument used to prove Theorem
\ref{thm:aftresult2} applies \emph{mutatis mutandis}.  Note that the cover $(E^\infty_f,{\mathcal L})$
is called the left Krieger cover in \cite{c,cm} and the Perron-Frobenius cover in \cite{sam}.  It is well
known that this cover is isomorphic to the left Krieger cover for two-sided sofic shifts.

\item  It is worth observing that for a two-sided sofic shift the cover $( E^\infty_f, {\mathcal L} )$
and the past set cover $( E_P , \mathcal{L}_P )$  are not necessarily isomorphic.
For instance, the following example was given in \cite[Section 4]{cm}:
Let $\textsf{Z}$ be the sofic shift
over the alphabet $\{ 1 , 2, 3 , 4 \}$ in which the forbidden blocks are
\[
\{ 1 2^k 1 , 3 2^k 1 2, 3 2^k 13 , 4 2^k 14 : k \ge 0 \}.
\]

\noindent The past set cover $( E_P , \mathcal{L}_P )$ of  $\textsf{Z}$ has five vertices while
the cover $(E^\infty_f,{\mathcal L})$ of $\textsf{Z}$ has four vertices.  The past set
cover contains $(E^\infty_f,{\mathcal L})$ as a subgraph.  Moreover, $(E_P,\mathcal{L}_P)$ is reducible,
but $(E^\infty_f,{\mathcal L})$ is irreducible.

\item  Further structural relationships between the five left-resolving covers of a sofic shift that we have described
in this paper are given in \cite{rune}.

\item  By recent results of Pask and Sunkara (see \cite{sunk}) the past set cover of a sofic $\beta$-shift is
irreducible. It follows by Theorem \ref{thm:aftresult2} that sofic $\beta$-shifts are not AFT.  This provides
a different proof of \cite[Proposition 2.7 c)]{thoms4}.  See also \cite{sam}.
\end{enumerate}
\end{rmks}


\begin{thebibliography}{llllll}
\bibitem{ash} Ashley, J. {\em Sliding block codes between constrained systems}, IEEE Trans. Inform. Theory {\bf 39} (1993), 1303--1309.
\bibitem{bprsz} Bates, T., Pask, D., Raeburn, I. and  Szyma{\' n}ski, W.
{\em The $C^*$-algebras of row-finite graphs}, New York J. Math. {\bf 6} (2000), 307--324.
\bibitem{bp2} Bates, T. and Pask, D. {\em The $C^*$-algebras of labelled graphs}, J. Operator Theory {\bf 57} (2007), 207--226.
\bibitem{b-m} Bertrand-Mathis, A. {\em D\'{e}veloppement en base $\theta$, r\'{e}partition modulo un de la suite $(x \theta^n )$, $n \ge 0$
langages cod\'{e}s et $\theta$-shift,} Bull. de la S.M.F. {\bf 114} (1986), 271-323.
\bibitem{bfp} B\'{e}al, M-P., Fiorenzi, F. and Perrin, D. {\em A hierarchy of shift equivalence sofic shifts}, Theoretical Computer Science, {\bf 345} (2005), 190--205.
\bibitem{blan} Blanchard, F. {\em Codes engendrant certains syst\`{e}mes sofiques}, Theoretical Computer Science, {\bf 68} (1989), 253--265.
\bibitem{bh} Boyle, M. and Huang, D.  {\em Poset block equivalence of integral matrices},  Trans. Amer. Math. Soc.  {\bf 355}  (2003), 3861--3886.
\bibitem{bkm} Boyle, M., Kitchens, B and Marcus, B. {\em A note on minimal covers for sofic systems}, Proc. Amer. Math. Soc. {\bf 95} (1985),  403--411.
\bibitem{bf} Bowen, R. and Franks, J. {\em Homology for zero-dimensional nonwandering sets},  Ann. Math. {\bf 106}  (1977), 73--92.
\bibitem{c} Carlsen, T. {\em  On $C^*$-algebras associated with sofic
shifts}, J.\ Operator Theory {\bf 49} (2003), 203--212.
\bibitem{ce1} Carlsen, T. and Eilers, S.  {\em Matsumoto $K$-groups associated to certain shift spaces},  Doc. Math.
{\bf 9}  (2004), 639--671.
\bibitem{ce2} Carlsen, T. and Eilers, S. {\em Augmenting dimension group invariants for substitution dynamics},  Ergod. Th. \& Dynam. Sys.  {\bf 24}  (2004), 1015--1039.
\bibitem{cm}  Carlsen, T. and Matsumoto, K. {\em Some remarks on the $C^*$-algebras associated with subshifts}, Math.\ Scand.\ {\bf
95} (2004), 145--160.
\bibitem{ffj} Fiebig, D., Fiebig, U-R. and Jonoska, N. {\em Multiplicities of covers for sofic shifts}, Theoretical Computer Science {\bf 262} (2001), 349--375.
\bibitem{fisc}  Fischer, R. {\em Sofic systems and graphs}, Monats. Math. {\bf 80} (1975), 179--186.
\bibitem{fo} Fujiwara, M. and Osikawa, M. {\em Sofic systems and flow equivalence}, Math. Rep. {\bf 16} (1987),
17--27.
\bibitem{rune} Johansen, R.  {\em On the structure of covers of irreducible sofic shifts}, in preparation.
\bibitem{jm} Jonoska, N. and Marcus, B. {\em Minimal presentations for irreducible sofic shifts}, IEEE Trans. Inform. Theory {\bf 40} (1994), 1818--1825.
\bibitem{kmw}  Katayama, Y., Matsumoto, K. and Watatani, Y. {\em Simple $C^*$-algebras arising from the $\beta$-expansion of real numbers,}
Ergod.\ Th.\ \& Dynam.\ Sys.\ {\bf 18} (1998), 937-962.
\bibitem{kit} Kitchens, B. ``Symbolic Dynamics:  One-sided, Two-sided and Countable State Markov Shifts", Springer, New York, 1998.
\bibitem{kri} Krieger, W. {\em On sofic systems I},  Israel J. Math., {\bf 48} (1984), 305--330.
\bibitem{kri1} Krieger, W. {\em On sofic systems II},  Israel J. Math.  {\bf 60}  (1987), 167--176.
\bibitem{kprr} Kumjian, A., Pask, D., Raeburn, I. and Renault, J. {\em Graphs, groupoids and Cuntz-Krieger algebras},
   J. Funct. Anal. {\bf 144} (1997), 505--541.
\bibitem{lm} Lind, D. and Marcus, B. ``An Introduction to Symbolic Dynamics and Coding",  Cambridge University Press, Cambridge, 1995.
\bibitem{marc} Marcus, B. {\em Sofic systems and encoding data}, IEEE Trans. Inform. Theory {\bf 31} (1985),  366--377.
\bibitem{m} Matsumoto, K. {\em On $C^*$-algebras associated with subshifts}, Internat.\ J.\ Math.\ {\bf 8}, (1997), 357-374.
\bibitem{m5} Matsumoto, K. {\em  $C^*$-algebras associated with presentations
of subshifts}, Documenta\ Math.\ {\bf 7}, (2002), 1-30.
\bibitem{m6} Matsumoto, K. {\em Bowen-Franks groups as an invariant for flow equivalence of subshifts},  Ergod. Th.
\& Dynam. Sys.  {\bf 21}  (2001),  1831--1842.
\bibitem{nasu} Nasu, M. {\em Topological conjugacy for sofic systems}, Ergod. Th. \& Dynam. Sys. {\bf 6} (1986),
265--280.
\bibitem{parsu} Parry, W. and Sullivan, D. {\em A topological invariant of flows on $1$-dimensional spaces},  Topology  {\bf 14}  (1975), 297--299.
\bibitem{sunk} Pask, D. and Sunkara, V. {\em Irreducible representations of certain $\beta$-shifts}, in preparation.
\bibitem{sam} Samuel, J. ``$C^*$-algebras of sofic shifts", PhD.\ Thesis Univ.\ Victoria, Canada (1998).
\bibitem{sha} Shannon, C. E., and Weaver, W. ``The Mathematical Theory of Communication", The University of Illinois Press, Urbana, Ill., 1949.
\bibitem{sen} Seneta, E. ``Non-negative matrices: An introduction to theory and applications",
\newblock Publ. George Allen \& Unwin Ltd., second edition (1973).
\bibitem{thoms3} Thomsen, K. {\em On the structure of a sofic shift space},  Trans. Amer. Math. Soc. {\bf 356} (2004), 3557--3619.
\bibitem{thoms4} Thomsen, K. {\em On the structure of beta shifts},  Algebraic and topological dynamics,  321--332, Contemp. Math., {\bf 385}, Amer. Math. Soc., Providence, RI, 2005.
\bibitem{trow} Trow, P. {\em Determining presentations of sofic shifts},  Theoretical Computer Science, {\bf 259} (2001), 199--216.
\bibitem{wei}  Weiss, B. {\em Subshifts of finite type and sofic systems}, Monats. Math. {\bf 77} (1973), 462--474.
\bibitem{will} Williams, S. {\em Covers of non-almost-finite type sofic systems}, Proc. Amer. Math. Soc. {\bf 104}
(1988),  245--252.
\end{thebibliography}
\end{document}